\documentclass[14pt,a4paper,fleqn]{article}
\usepackage{amsmath,amssymb,amsfonts}
\newtheorem{theorem}{Theorem}

\numberwithin{equation}{section}
\numberwithin{lemma}{section}
\numberwithin{theorem}{section}
\numberwithin{corollary}{section}

\allowdisplaybreaks

\begin{document}
\title{Infinite summation formulas of Srivastava's general triple hypergeometric function}
\author{Vivek Sahai
\\ 
Department of Mathematics and Astronomy\\ Lucknow University \\Lucknow 226007, India 
\\
sahai\_vivek@hotmail.com\\[10pt]
Ashish Verma\footnote{Corresponding author}
\\ 
Department of Mathematics\\ Prof. Rajendra Singh (Rajju Bhaiya)\\ Institute of Physical Sciences for Study and Research \\  V. B. S. Purvanchal University, Jaunpur  (U.P)- 222003, India\\
vashish.lu@gmail.com}
\maketitle
\begin{abstract}In this paper we derive the  infinite summation  formulas of Srivastava's general triple hypergeometric function. Certain particular cases leading to  infinite summation formulas for fourteen Lauricella and three Srivastava\'s triple hypergeometric functions are also presented.\\[12pt]
Keywords:  Srivastava's general triple hypergeometric function, Lauricella function, Srivastava’s triple hypergeometric  function.
\end{abstract}
\section{Introduction}
Recently, the authors have studied the Srivastava's general triple hypergeometric function \cite {VS, VS1} from the view point of recursion formulas and $q$-derivative with respect to parameters. Further, in  \cite{VS2}, we have enumerated finite summation formulas of Srivastava's general triple hypergeometric function.  In the present paper, we obtain the infinite summation formulas satisfied by this function. The particular cases will lead to results involving fourteen Lauricella function and three Srivastava's triple hypergeometric function.\\

Earlier,  Wang \cite{XW} obtained infinite summation  formulas of   double hypergeometric functions. In \cite{XW1}, finite summation formulas of double hypergeometric functions are derived.

The Srivastava's general triple hypergeometric function \cite{S, SK,SM} is defined by 
\begin{eqnarray}
F^{(3)}[x, y, z]&=&F^{(3)}\left [^{(a)\,\,::\,\,\,(b); \,\,\,(b')\,\,; \,\,\,(b'')\,\,\,:\,\,\,(c)\,\,; \,\,\,(c')\,\,; \,\,\,(c'') \,\,;}_{(e)\,\,::\,\,\,(g)\,; \,\,(g')\,\,;\,\,\, (g'')\,\,:\,\,\,(h);\,\,\,\,(h');\,\,\,(h'')\,\,\,;} \,\,x_1, x_2, x_3\right ]\notag\\
&=&\sum_{m_1, m_2, m_3= 0}^{\infty}\wedge(m_1,\,m_2,\,m_3)\frac{{x_1}^{m_1}}{{m_1}{!}}\frac{{x_2}^{m_2}}{{m_2}{!}}\frac{{x_3}^{m_3}}{{m_3}{!}}\label{1eq2}
\end{eqnarray}\\
where  
\begin{eqnarray}
&&\wedge(m_1,\,m_2,\,m_3)\notag\\
&=& \frac{\displaystyle\prod_{j=1}^{A}(a_j)_{m_1+m_2+m_3} \prod_{j=1}^{B}(b_j)_{m_1+m_2}\prod_{j=1}^{B'}(b_{j}')_{m_2+m_3}\prod_{j=1}^{B''}(b_{j}'')_{m_3+m_1}\prod_{j=1}^{C}(c_{j})_{m_1}\prod_{j=1}^{C'}(c_{j}')_{m_2}\prod_{j=1}^{C''}(c_{j}'')_{m_3}}{\displaystyle\prod_{j=1}^{E}(e_j)_{m_1+m_2+m_3} \prod_{j=1}^{G}(g_j)_{m_1+m_2}\prod_{j=1}^{G'}(g_{j}')_{m_2+m_3}\prod_{j=1}^{G''}(g_{j}'')_{m_3+m_1}\prod_{j=1}^{H}(h_{j})_{m_1}\prod_{j=1}^{H'}(h_{j}')_{m_2}\prod_{j=1}^{H''}(h_{j}'')_{m_3}}\notag\\ \label{1eq3}
\end{eqnarray}
and $(a)$ abbreviates the array of $A$ parameters $a_1, a_2,\ldots, a_{A}$, etc.  The region of convergence of the general triple hypergeometric series (\ref{1eq2}) is given in the literature \cite{SM}. \\

In this paper following abbreviated notations are used. For example, write 
\begin{align}
&(a+k):=a_1+k,\dots,a_A+k,\notag\\
&(a^j):=a_1,\dots,a_{j-1}, a_{j+1},\dots,a_{A},\,\,\notag\\
&(a^j+k):=a_{1}+k,\dots,a_{j-1}+k, a_{j+1}+k,\dots,a_{A}+k,\,\,\, j=1,\dots,A.\notag\end{align}
Also, we use the notations
 \begin{align}
  \qquad [a]_k:=\prod_{i=1}^{A}\,(a_i)_{k}
 ,\qquad [a^j]_k:=\prod_{i=1,i\neq j}^{A}\,(a_i )_k ,\,\, etc. \notag
\end{align}
where $k$ is a non-negative integer and  $(a_i)_k$ is the Pochhammer symbol \cite{R}.\\

The subject matter is divided in two sections. The infinite summation formulas for non-terminating Srivastava's general triple hypergeometric function are presented in Section~2 and the results for  the terminating Srivastava's general triple hypergeometric function are discussed in Section~3. 
\section{Infinite summation formulas with non-terminating Srivastava's general triple hypergeometric function}
In this section, we establish infinite summation formulas of non-terminating Srivastava's general triple hypergeometric function by using the well known binomial summation theorem.
\begin{theorem}\label{t1}
The following infinite summation formulas of Srivastava's general triple hypergeometric function hold true:
\begin{align}
&\sum_{k=0}^{\infty}\frac{(a_{i})_{k}}{k!}\,t^{k}F^{(3)}\left [^{\,a_{i}+k, \,(a^{i})\,::\,(b)\,;\,(b')\,; \,(b''): \,(c)\,; (c')\,; \,(c'') \,;}_{\,\,\,\,\,\,\,\,\,(e)\,\,\,\,\,\,\,::\,(g)\,; \,(g')\,;\, (g'')\,:\,(h)\,;\,(h');\,(h'')\,;} \,\,x_1, x_2, x_3\right ]\notag\\
&=\, (1-t)^{-a_{i}}\,F^{(3)}\left [^{(a)::\,(b)\,; \,(b')\,; \,(b''):\,(c); \,(c')\,; (c'') ;}_{(e)\,::\,(g)\,;\,(g')\,;\, (g''):\,(h)\,;\,(h');\,(h'')\,;} \,\frac{x_1}{1-t}, \frac{x_2}{1-t}, \frac{x_3}{1-t}\right ],\label{1s2}
\end{align}
where $i=1,\dots,A$;
\begin{align}
&\sum_{k=0}^{\infty}\frac{(b_{i})_{k}}{k!}\,t^{k}F^{(3)}\left [^{(a)\,\,::\,b_{i}+k, \,(b^{i})\,;\,(b')\,; \,(b''): \,(c)\,; (c')\,; \,(c'') \,;}_{(e)\,\,::\,\,\,\,\,\,\,\,\,\,(g)\,\,\,\,\,\,; \,(g')\,;\, (g'')\,:\,(h)\,;\,(h');\,(h'')\,;} \,\,x_1, x_2, x_3\right ]\notag\\
&=\, (1-t)^{-b_{i}}\,F^{(3)}\left [^{(a)::\,(b)\,; \,(b')\,; \,(b''):\,(c); \,(c')\,; (c'') ;}_{(e)\,::\,(g)\,;\,(g')\,;\, (g''):\,(h)\,;\,(h');\,(h'')\,;} \,\frac{x_1}{1-t}, \frac{x_2}{1-t}, x_3\right ],\label{t1s2}
\end{align}
where $i=1,\dots,B$;
\begin{align}
&\sum_{k=0}^{\infty}\frac{(c_{i})_{k}}{k!}\,t^{k}F^{(3)}\left [^{(a)\,\,::\,(b); \,(b'); \,(b''):\,c_{i}+k, (c^{i})\,; (c')\,; \,(c'') \,;}_{(e)\,\,::\,(g)\,; \,(g')\,;\, (g'')\,:\,\,\,\,\,(h)\,\,\,\,\,\,\,\,;\,(h');\,(h'')\,;} \,\,x_1, x_2, x_3\right ]\notag\\
&=\, (1-t)^{-c_{i}}\,F^{(3)}\left [^{(a)::\,(b); (b')\,; \,(b''):\,(c); \,(c')\,; (c'') ;}_{(e)\,::\,(g)\,;\,(g')\,;\, (g''):\,(h)\,;\,(h');\,(h'')\,;} \,\frac{x_1}{1-t}, x_2, x_3\right ],\label{t2s2}
\end{align}
where $i=1,\dots,C$.
\end{theorem}
{\bf Proof:}  Applying the definition of Srivastava's general triple hypergeometric function $F^{\left(3\right)}{[x, y, z]}$ and  transformation \begin{align*}(a_{i})_{k} \,(a_{i}+k)_{m_1+m_2+m_3}= (a_{i})_{m_1+m_2+m_3}\, (a_{i}+m_1+m_2+m_3)_{k},\end{align*} the L.H.S  of (\ref{1s2}) can be expressed as  
\begin{align*}
\sum_{m_1, m_2, m_3=0}^{\infty}\wedge(m_1, m_2, m_3)\frac{x_{1}^{m_1} x_{2}^{m_2} x_{3}^{m_3}}{m_{1}! m_2! m_3!}\, _{1}F_{0}\left[^{a_{i}+m_1+m_2+m_3}_{\,\,\,\,\quad-\,\,\,\,}; t\right],
\end{align*}
where 
 \cite{SM}
\begin{align}
_{1}F_{0}\left[^{a}_{\,-\,}; t\right]= (1-t)^{-a}.\label{N}
\end{align}
After some simplification, we get the right side of (\ref{1s2}). This completes the proof of (\ref{1s2}). The remaining identities (\ref{t1s2}) and (\ref{t2s2}) are proved in a similar manner.\\

\begin{theorem}
The following infinite summation formulas of Srivastava's general triple hypergeometric function hold true:
\begin{align}
F^{(3)}\left [x_1+t, x_2, x_3\right ]=&\sum_{k=0}^{\infty}\frac{[a]_{k}[b]_{k}[b'']_{k}[c]_{k}}{[e]_{k}[g]_{k}[g'']_{k}[h]_{k} k!}\,\,t^{k}\notag\\
\qquad\qquad&\times F^{(3)}\left [^{(a+k)\,\,::\,(b+k)\,;\,(b')\,; \,(b''+k)\,:\,(c+k)\,; \,(c')\,; \,(c'')\, ;}_{\,(e+k)\,\,::\,(g+k)\,; \,(g')\,;\,(g''+k)\,:\,(h+k)\,;\,(h')\,;\,(h'')\,;} \,\,x_1, x_2, x_3\right ];\label{s3}
\end{align}
\begin{align}
F^{(3)}\left [x_1, x_2+t, x_3\right ]=&\sum_{k=0}^{\infty}\frac{[a]_{k}[b]_{k}[b']_{k}[c']_{k}}{[e]_{k}[g]_{k}[g']_{k}[h']_{k} k!}\,\,t^{k}\notag\\
\qquad\qquad&\times F^{(3)}\left [^{(a+k)\,\,::\,(b+k)\,;\,(b'+k)\,; \,(b'')\,:\,(c)\,; \,(c'+k)\,; \,(c'')\, ;}_{\,(e+k)\,\,::\,(g+k)\,; \,(g'+k)\,;\,(g'')\,:\,(h)\,;\,(h'+k)\,;\,(h'')\,;} \,\,x_1, x_2, x_3\right ];\label{t1s3}
\end{align}
\begin{align}
F^{(3)}\left [x_1, x_2, x_3+t\right ]=&\sum_{k=0}^{\infty}\frac{[a]_{k}[b']_{k}[b'']_{k}[c'']_{k}}{[e]_{k}[g']_{k}[g'']_{k}[h'']_{k} k!}\,\,t^{k}\notag\\
\qquad\qquad&\times F^{(3)}\left [^{(a+k)\,\,::\,(b)\,;\,(b'+k)\,; \,(b''+k)\,:\,(c)\,; \,(c')\,; \,(c''+k)\, ;}_{\,(e+k)\,\,::\,(g)\,; \,(g'+k)\,;\,(g''+k)\,:\,(h)\,;\,(h')\,;\,(h''+k)\,;} \,\,x_1, x_2, x_3\right ].\label{t2s3}
\end{align}
\end{theorem}
{\bf Proof:} From the definition $F^{\left(3\right)}{[x_1, x_2, x_3]}$ and the transformation $(a)_{k}\,(a+k)_{m_1}= (a)_{k+m_1}$, the right side of (\ref{s3}) can be written as 
\begin{align*}
\sum_{k, m_1, m_2, m_3=0}^{\infty}\wedge(m_1+k, m_2, m_3)\frac{x_{1}^{m_1} x_{2}^{m_2} x_{3}^{m_3}}{m_{1}! m_{2}! m_{3}! k!}\, t^{k}.
\end{align*}
Replacing $m_{1}+k\rightarrow l$ in the above result and after some simplification, we get
\begin{align*}
&\sum_{l, m_{2}, m_{3}=0}^{\infty}\wedge(l, m_{2}, m_{3})\frac{ x_{2}^{m_{2}} x_{3}^{m_{3}}}{l! m_{2}! m_{3}!}\, \sum_{k=0}^{l}{l\choose k}x_{1}^{l-k} t^{k}\notag\\
&=\sum_{l, m_{2}, m_{3}=0}^{\infty}\wedge(l, m_{2}, m_{3})\frac{(x_{1}+t)^{l} x_{2}^{m_{2}} x_{3}^{m_{3}}}{l! m_{2}! m_{3}!}.
\end{align*}
where, we have applied the special case of binomial theorem (\ref{N}) as 
\begin{align}
\sum_{k=0}^{l}{l\choose k} x_{1}^{k} x_{2}^{l-k}= (x_1+x_2)^{l}
\end{align}
in the inner summation. This completes the proof of (\ref{s3}).  Identities (\ref{t1s3}) to  (\ref{t2s3}) can be proved in an analogous manner.\\

\begin{theorem}
The following infinite summation formulas of Srivastava's general triple hypergeometric function hold true:
\begin{align}
&\sum_{k=0}^{\infty}\frac{[a^{i}]_{k}[b]_{k}[b'']_{k}[c]_{k}(r)_{k}}{[e]_{k}[g]_{k}[g'']_{k}[h]_{k}k!}\, x_{1}^{k}\,F^{(3)}\left [^{\, a_{i},\,(a^{i}+k)\,\,::\,(b+k)\,; \,(b')\,; \,(b''+k)\,:\,(c+k)\,; \,(c')\,; \,(c'') \,;}_{\,\,\,\,\,\,\,\,(e+k)\,\,\,\,\,\,::\,(g+k)\,; \,(g')\,;\, (g''+k)\,:\,(h+k)\,;\,(h');\,(h'');} \,x_1, x_2, x_3\right ]\notag\\
&=\,F^{(3)}\left [^{\,a_{i}+r,\,\,(a^{i})\,::\,(b)\,;\,\,\,\,a_{i},(b')\,; \,(b'')\,: \,(c)\,\,\,\,\,\,; \,(c')\,; \,(c'') \,;}_{\,\,\,\,\,\,\,\,\,\,\,\,\,\,(e)\,\,\,\,\,\,::\,(g)\,; a_{i}+r, (g')\,;\, (g'')\,:\,\,(h)\,;\,(h');\,(h'')\,;} \,x_1, x_2, x_3\right ],\label{s4}
\end{align}
where $i=1,\dots,A;$
\begin{align}
&\sum_{k=0}^{\infty}\frac{[a]_{k}[b]_{k}[b'']_{k}[c^{i}]_{k}(r)_{k}}{[e]_{k}[g]_{k}[g'']_{k}[h]_{k}k!}\, x_{1}^{k}\,F^{(3)}\left [^{(a+k)\,\,::\,(b+k)\,; \,(b')\,; \,(b''+k)\,:\,c_{i}, (c^{i}+k)\,; \,(c')\,; \,(c'') \,;}_{\,(e+k)\,::\,(g+k)\,; \,(g')\,;\, (g''+k)\,:\,\,\,\,\,(h+k)\,\,\,\,;\,(h');\, (h'')\,;} \,\,x_1, x_2, x_3\right ]\notag\\
&=\,F^{(3)}\left [^{(a)\,\,::\,(b)\,;\,(b')\,; \,(b'')\,:\,c_{i}+r,\, (c^{i})\,; \,(c')\,; \,(c'') \,;}_{(e)\,\,::\,(g)\,; \,(g')\,;\, (g'')\,:\,\,\,\,\,\,\,(h)\,\,\,\,\,\,\,\,\,\,;\,(h');\,(h'')\,;} \,x_1, x_2, x_3\right ],\label{t1s4}
\end{align}
where $i=1,\dots,C.$
\end{theorem}
{\bf Proof:} We first prove identity (\ref{s4}). 
 From the definition of $F^{\left(3\right)}{[x_1, x_2, x_3]}$-series and  the transformation $(a)_{k} (a+k)_{m_1}= (a)_{k+m_1}$, the left side of (\ref{s4}) can be expressed as
\begin{align*}
&\sum_{k,m_1, m_2, m_3=0}^{\infty}\frac{\wedge(m_1+k, m_2, m_3)(r)_{k}}{(a_{i}+m_1+m_2+m_3)_{k}\,k!}\frac{x_{1}^{m_1+k} x_{2}^{m_2} x_{3}^{m_3}}{m_{1}! m_{2}! m_{3}!}\notag\\
&=\sum_{l, m_{2}, m_{3}=0}^{\infty}\wedge(l, m_{2}, m_3)\frac{x_{1}^{l} x_{2}^{m_{2}} x_{3}^{m_{3}}}{l! m_{2}! m_{3}!}\,_{2}F_{1}\left[^{\,\,\,\,-l,\,\, r\,\,\,}_{1-a_{i}-l-m_{2}-m_{3}}; 1\right],
\end{align*}
where, we have performed the replacement $m_1+k\rightarrow l$ in the second equation. Obviously the inner summation $_{2}F_{1}$ can be evaluated by the Vandermonde's theorem \cite{R}
\begin{align}
_{2}F_{1}\left[^{-n, a\,}_{\,\,\,\,c\,\,\,\,}; 1\right]=\frac{(c-a)_{n}}{(c)_{n}}\label{V}
\end{align}
After some simplification, we get right side of (\ref{s4}). The identity (\ref{t1s4}) can be proved in an analogous manner. 
\begin{theorem}\label{t2}
The following infinite summation formulas of Srivastava's general triple hypergeometric function hold true:
\begin{align}
&\sum_{k=0}^{\infty}\frac{[a^{i}]_{k}[b]_{k}[b'']_{k}[c]_{k}(d)_{k}}{[e]_{k}[g]_{k}[g'']_{k}[h]_{k}\,k!}\,(-x_1)^{k}\,F^{(3)}\left [^{\,(a+k)\,::\,(b+k)\,; \,(b')\,; \,(b''+k)\,:\,(c+k)\,; \,(c')\,; \,(c'')\, ;}_{ \,(e+k)\,::\,(g+k)\,; \,(g')\,;\, (g''+k)\,:\,(h+k)\,;\,(h')\,;\,(h'');} {x_1}, {x_2}, {x_3}\right ]\notag\\
&=\,F^{(3)}\left [^{\,(a)\,::\,(b)\,; \,(b')\,; \,(b'')\,:\,\,a_{i}-d,\,(c)\,; \,(c')\,;\,(c'')\,;}_{\,\,(e)\,\,::\,(g)\,; \,(g')\,;\, (g'')\,:\,\,\,a_i, (h)\,\,\,\,\,;\,(h');\,(h'')\,;} x_1, x_2, x_3\right ],\label{1s6}
\end{align}
where $i=1,\dots, A;$
\begin{align}
&\sum_{k=0}^{\infty}\frac{[a]_{k}[b]_{k}[b'']_{k}[c^{i}]_{k}(d)_{k}}{[e]_{k}[g]_{k}[g'']_{k}[h]_{k}\,k!}\,(-x_1)^{k}\,F^{(3)}\left [^{(a+k)\,::\,(b+k)\,;\, (b')\,; \,(b''+k)\,:\,(c+k)\,; \,(c')\,;\,(c'')\,;}_{(e+k)\,::\,(g+k)\,;\, (g')\,;   \,(g''+k)\,:\,(h+k)\,\,;\,(h');\,(h'')\,;} x_1, x_2, x_3 \right ]\notag\\
&=\,F^{(3)}\left [^{(a)\,::\,(b); (b')\,;\,(b'')\,:\,c_{i}-d,\,(c^{i})\,; \,(c')\,; \,(c'')\,;}_{(e)\,::\,(g)\,; \,(g')\,; \,(g'')\,:\,\,\,\,\,\,\,\,\,(h)\,\,\,\,\,\,;\,(h')\,;\,(h'')\,;} x_1, x_2, x_3\right ],\label{m2s6}
\end{align}
where $i=1,\dots, C.$
\end{theorem}
{\bf Proof:} We give the proof of (\ref{1s6}). 
This theorem can be proved by Vandermonde's theorem (\ref{V}). Using the definition of $F^{(3)}[x_1, x_2, x_3]$ and  the transformation $(a)_{k} (a+k)_{m_{1}}= (a)_{k+m_{1}}$, then applying replacement $m_1+k\rightarrow l$, and making some simplification, the left side of (\ref{1s6}) can be expressed as
\begin{align*}
\sum_{l, m_2, m_3=0}^{\infty}\frac{\wedge(l, m_2, m_3)}{l!\, m_2!\, m_3!} x_{1}^{l} x_{2}^{m_{2}} x_{3}^{m_{3}}\, _{2}F_{1}\left[^{\,-l, d\,\,}_{\,\,\,\,\,\,a_{i}\,\,}; 1\right].
\end{align*}
Calculating inner summation $_{2}F_{1}(1)$ by Vandermonde's theorem (\ref{V}), we get right side of (\ref{1s6}). Transformation (\ref{m2s6})  can be proved in an analogous manner. We omit these details.

\begin{theorem}The following infinite summation formulas of Srivastava's general triple hypergeometric function hold true:
\begin{align}\
&\sum_{k=0}^{\infty}\frac{[a]_{k}[b]_{k}[b'']_{k}[c^{i}]_{k}(d)_{k}(r)_{k}}{[e]_{k}[g]_{k}[g'']_{k}[h]_{k}(d+r+c_{i})_{k}\,k!}\,x_{1}^{k}\,F^{(3)}\left [^{(a+k)\,::\,(b+k)\,;\, (b')\,; \,(b''+k)\,:\,c_{i},\,(c^{i}+k)\,; \,(c')\,;\,(c'')\,;}_{(e+k)\,::\,(g+k)\,;\, (g')\,;   \,(g''+k)\,:\,\,\,\,\,(h+k)\,\,\,\,\,\,;\,(h');\,(h'')\,;} x_1, x_2, x_3 \right ]\notag\\
&=\,F^{(3)}\left [^{(a)\,::\,(b); (b')\,;\,(b'')\,:\,c_{i}+r,\,c_{i}+d,\,(c^{i})\,; \,(c')\,; \,(c'')\,;}_{(e)\,::\,(g)\,; \,(g')\,; \,(g'')\,:\,\,\,c_{i}+r+d,\,\,(h)\,\,\,;\,(h')\,;\,(h'')\,;} x_1, x_2, x_3\right ],\label{l2s6}
\end{align}
where $i=1,\dots, C.$
\end{theorem}
{\bf Proof:} The left side of (\ref{l2s6}) can be expressed as
\begin{align*}
&\sum_{k,m_1, m_2, m_3=0}^{\infty}\frac{\wedge(m_1+k, m_2, m_3)(r)_{k}(d)_{k}}{m_{1}! m_{2}! m_{3}! k!}\frac{(c_i)_{m_{1}}\,x_{1}^{m_{1}+k} x_{2}^{m_2} x_{3}^{m_{3}}}{(c_i)_{m_1+k}(d+r+c_i)_{k} } \notag\\
&= \sum_{l, m_2, m_3=0}^{\infty}\frac{\wedge(l, m_2, m_3)}{l! m_2! m_3!}x_{1}^{l} x_{2}^{m_{2}} x_{3}^{m_3}\, _{3}F_{2}\left[^{\,\,\,-l,\,\, d,\,\,\,r\,\,}_{d+r+c_i, 1-c_i-l}; 1\right],
\end{align*}
where, Replaced $m_1+k\rightarrow l$ in above equation and using the well known Saalschutz  formula \cite{L}
\begin{align}
_{3}F_{2}\left[^{\,\,\,-n,\,\, a,\,\,b\,\,}_{c, 1+a+b-c-n}; 1\right]=\frac{(c-a)_{n}(c-b)_{n}}{(c)_{n}(c-a-b)_{n}}\label{S}
\end{align}
and after some simplification, we get right side of (\ref{l2s6}). This completes the proof of the theorem.
\begin{theorem} The following infinite summation formulas of Srivastava's general triple hypergeometric function hold true:
\begin{align}
&\sum_{k=0}^{\infty}\frac{[a^{i}]_{k}[b]_{k}[b'']_{k}[c]_{k}(d)_{k}(1+\frac{d}{2})_{k}}{[e]_{k}[g]_{k}[g'']_{k}[h]_{k}(\frac{d}{2})_{k}\,k!}\,(-x_1)^{k}\,F^{(3)}\left [^{(a+k)\,::\,(b+k)\,;\, (b')\,; \,(b''+k)\,:\,(c+k)\,; \,(c')\,;\,(c'')\,;}_{(e+k)\,::\,(g+k)\,;\, (g')\,;   \,(g''+k)\,:\,(h+k)\,;\,(h');\,(h'')\,;}x_1, x_2, x_3 \right ]\notag\\
&=\,F^{(3)}\left [^{(a)\,::\,(b); (b')\,;\,(b'')\,:\,2+d-a_{i},\,a_{i}-d-1,\,(c)\,; \,(c')\,; \,(c'')\,;}_{(e)\,::\,(g)\,; \,(g')\,; \,(g'')\,:\,\,\,\,\,\,\,1+d-a_{i},\,\,\,a_{i},\,\,\,\,(h)\,\,;\,(h')\,;\,(h'')\,;} x_1, x_2, x_3\right ],\label{n2s6}
\end{align}
where $i=1,\dots, A;$
\begin{align}\
&\sum_{k=0}^{\infty}\frac{[a]_{k}[b]_{k}[b'']_{k}[c^{i}]_{k}(d)_{k}(1+\frac{d}{2})_{k}}{[e]_{k}[g]_{k}[g'']_{k}[h]_{k}(\frac{d}{2})_{k}\,k!}\,(-x_1)^{k}\,\,F^{(3)}\left [^{(a+k)\,::\,(b+k)\,;\, (b')\,; \,(b''+k)\,:\,(c+k)\,; \,(c')\,;\,(c'')\,;}_{(e+k)\,::\,(g+k)\,;\, (g')\,;   \,(g''+k)\,:\,(h+k)\,;\,(h');\,(h'')\,;} x_1, x_2, x_3 \right ]\notag\\
&=\,F^{(3)}\left [^{(a)\,::\,(b); (b')\,;\,(b'')\,:\,2+d-c_{i},\,c_{i}-d-1,\,(c^{i})\,; \,(c')\,; \,(c'')\,;}_{(e)\,::\,(g)\,; \,(g')\,; \,(g'')\,:\,\,\,\,\,\,\,1+d-c_{i},\,\,\,\,\,(h)\,\,\,\,\,\,\,\,\,\,\,;\,(h')\,;\,(h'')\,;} x_1, x_2, x_3\right ],\label{t2s6}
\end{align}
where $i=1,\dots, C.$
\end{theorem}
{\bf Proof:} The left side of (\ref{n2s6}) can be expressed as 
\begin{align*}
&\sum_{k, m_1, m_2, m_3=0}^{\infty}\frac{\wedge(m_1+k, m_2, m_3)(d)_{k}(1+\frac{d}{2})}{m_1!\,m_2!\,m_3! k!(\frac{d}{2})_{k}(a_{i})_{k}(-1)^{k}}\, x_{1}^{m_1+k} x_{2}^{m_2} x_{3}^{m_3}\notag\\
&=\sum_{l, m_{2}, m_{3}=0}^{\infty}\frac{\wedge(l,m_{2}, m_{3})}{l! m_{2}! m_{3}!}\, x_{1}^l x_{2}^{m_2} x_{3}^{m_3}\,_{3}F_{2}\left[^{\,\,\,-l,\,\, d,\,1+\frac{d}{2}\,\,}_{\,\,\,\,\,\,\,\frac{d}{2}, \,\,\,\,\,\,a_{i}}; 1\right],
\end{align*}
where replace $m_1+k\rightarrow l$ in the last equation. Applying the  nearly poised summation formula \cite{L}:
\begin{align}
_{3}F_{2}\left[^{\,\,\,-n,\,\, a,\,1+\frac{a}{2}\,\,}_{\,\,\,\,\frac{a}{2},\,\, b}; 1\right]=(b-a-1-n)\frac{(b-a)_{n-1}}{(b)_{n}}\label{np}
\end{align}
and simplifying, we get (\ref{n2s6}). The transformation (\ref{t2s6})  can be proved in an analogous manner. 
\begin{theorem}The following infinite summation formula of Srivastava's general triple hypergeometric function holds true:
\begin{align}
&\sum_{k=0}^{\infty}\frac{[a]_{k}[b]_{k}[b'']_{k}[c^{i}]_{k}(r)_{k}\,(\frac{-c_{i}}{2})_{k}}{[e]_{k}[g]_{k}[g'']_{k}[h]_{k}(1+r+\frac{c_{i}}{2})_{k}\,k!}\,x_{1}^{k}\,\,F^{(3)}\left [^{(a+k)\,::\,(b+k)\,;\, (b')\,; \,(b''+k)\,:\,c_{i},\,(c^{i}+k)\,; \,(c')\,;\,(c'')\,;}_{(e+k)\,::\,(g+k)\,;\, (g')\,;   \,(g''+k)\,:\,\,\,\,\,\,(h+k)\,\,\,\,\,;\,(h');\,(h'')\,;} x_1, x_2, x_3 \right ]\notag\\
&=\,F^{(3)}\left [^{(a)\,::\,(b); (b')\,;\,(b'')\,:\,c_{i}+r,\,\frac{c_{i}}{2},\,1+\frac{1}{2}(c_{i}+r),\,(c^{i})\,; \,(c')\,; \,(c'')\,;}_{(e)\,::\,(g)\,; \,(g')\,; \,(g'')\,:\,\,1+r+\frac{c_{i}}{2},\,\frac{1}{2}(c_{i}+r),\,(h)\,\,\,\,\,;\,(h')\,;\,(h'')\,;} x_1, x_2, x_3 \right ],\label{2s7}
\end{align}
where $i=1,\dots, C.$
\end{theorem}
{\bf Proof:} The left side of (\ref{2s7}) can be expressed as 
\begin{align*}
&\sum_{k, m_1, m_2, m_3=0}^{\infty}\frac{\wedge(m_1+k, m_2, m_3)(r)_{k}(\frac{-c_i}{2})_{k}}{m_{1}! m_{2}! m_{3}! k!(1+r+\frac{c_i}{2})_{k}}\frac{(c_i)_{m_1}}{(c_i)_{m_1+k}}\, x_{1}^{m_{1}+k} x_{2}^{m_2} x_{3}^{m_3}\,\notag\\
&=\sum_{l, m_2, m_3=0}^{\infty}\frac{\wedge(l, m_2, m_3)}{l! m_2! m_3!}\,x_{1}^{l} x_{2}^{m_2} x_{3}^{m_3} \,_{3}F_{2}\left[^{\,\,\,-l,\,\, r,\,\frac{-c_i}{2}\,\,}_{\,1+r+\frac{c_i}{2},\,\, 1-c_i-l}; 1\right].
\end{align*}
where, we have replaced $m_1+k\rightarrow l$ in the second equation. Evaluating the inner summation using the summation formula \cite{L} 
\begin{align}
_{3}F_{2}\left[^{\,\,\,-n,\,\, a,\,b\,\,}_{\,1+a-b,\,\, 1+2b-n}; 1\right]=
\frac{(a-2b)_n(1+\frac{1}{2}a-b)_n(-b)_n}{(1+a-b)_n (\frac{a}{2}-b)_n(-2b)_n}
\end{align}
and simplifying, we get the right side  of (\ref{2s7}). This completes the proof of this theorem.
\begin{theorem}The following infinite summation formula of Srivastava's general triple hypergeometric function holds true:
\begin{align}
&\sum_{k=0}^{\infty}\frac{[a]_{k}[b]_{k}[b'']_{k}[c^{i}]_{k}(d)_{k}\,(1+\frac{d}{2})_{k}(\frac{-c_{i}}{2})_{k}}{[e]_{k}[g]_{k}[g'']_{k}[h]_{k}(1+d+\frac{c_{i}}{2})_{k}(\frac{d}{2})_{k}\,k!}\,x_{1}^{k}\,\,F^{(3)}\left [^{(a+k)\,::\,(b+k)\,;\, (b')\,; \,(b''+k)\,:\,c_{i},\,(c^{i}+k)\,; \,(c')\,;\,(c'')\,;}_{(e+k)\,::\,(g+k)\,;\, (g')\,;   \,(g''+k)\,:\,\,\,\,\,\,(h+k)\,\,\,\,\,;\,(h');\,(h'')\,;} x_1, x_2, x_3 \right ]\notag\\
&=\,F^{(3)}\left [^{(a)\,::\,(b); (b')\,;\,(b'')\,:\,\frac{c_{i}}{2},\, c_{i}+d,\,(c^{i})\,; \,(c')\,; \,(c'')\,;}_{(e)\,::\,(g)\,; \,(g')\,; \,(g'')\,:\,1+d+\frac{c_{i}}{2},\,(h)\,;\,(h')\,;\,(h'');} x_1, x_2, x_3\right ],\label{t2s7}
\end{align}
where $i=1,\dots, C.$
\end{theorem}
{\bf Proof:} From the definition of $F^{(3)}[x_1, x_2, x_3]$-series and using the transformation $(a)_{k}(a+k)_{m_1}= (a)_{k+m_{1}}$, the left side of (\ref{t2s7}) can be simplified as 
\begin{align*}
&\sum_{k, m_1, m_2, m_3=0}^{\infty}\frac{\wedge(m_1+k, m_2, m_3)}{m_{1}! m_2! m_3! k!}\frac{(c_i)_{m_1}(d)_{k}\,(1+\frac{d}{2})_{k}(\frac{-c_{i}}{2})_{k}}{(c_i)_{m_{1}+k}(1+d+\frac{c_{i}}{2})_{k}(\frac{d}{2})_{k}}\, x_{1}^{m_{1}+k} x_{2}^{m_2} x_{3}^{m_{3}}\notag\\
&=\sum_{l, m_2, m_3=0}^{\infty}\frac{\wedge(l, m_2, m_3)}{l! m_2! m_3!}\, x_{1}^l x_{2}^{m_3} x_{3}^{m_3}\, _{4}F_{3}\left[^{\,-l,\,d,\,1+\frac{d}{2}, -\frac{1}{2} c_{i}\,\,}_{\,\frac{1}{2}d,\,1+d+\frac{1}{2}c_{i}\,, 1-c_{i}-l}; 1\right]
\end{align*} 
where we have replaced $m_1+k\rightarrow l$ in the second equation. Applying the following well known summation formula \cite{L} in the inner summation $_{4}F_{3}(1)$-series
\begin{align}
_{4}F_{3}\left[^{\,-n,\,a,\,1+\frac{a}{2}, b\,\,}_{\,\frac{1}{2}a,\,1+a-b\,, 1+2b-n}; 1\right]=\frac{(a-2b)_{n}(-b)_{n}}{(1+a-b)_{n}(-2b)_{n}},
\end{align}
 and simplifying, we get (\ref{t2s7}).

\section{The infinite summation formulas with terminating Srivastava's general triple hypergeometric function}
In this section, we present infinite summation formulas with terminating Srivastava's general triple hypergeometric function by using the binomial theorem (\ref{N}).
\begin{theorem}\label{t9}
The following infinite summation formula of Srivastava's general triple hypergeometric function holds true:
\begin{align}
&\sum_{k=0}^{\infty}\frac{(c_{i})_{k}}{k!}\,(-t)^{k}F^{(3)}\left [^{\,(a)\,::\,(b)\,;\,(b')\,; \,(b''): \,-k, (c^{i})\,; (c')\,; \,(c'') \,;}_{\,(e)\,::\,(g)\,; \,(g')\,;\, (g'')\,:\,\,\,(h)\,\,\,\,\,;\,(h');\,(h'')\,;} \frac{1+t}{t} x_1, x_2, x_3\right ]\notag\\
&=\, (1+t)^{-c_{i}}\,F^{(3)}\left [ x_1, x_2, x_3\right],\label{10s1}
\end{align}
where $i=1,\dots,C$.
\end{theorem}
{\bf Proof:} From the definition of $F^{(3)}[x_1, x_2, x_3]$, the left side of (\ref{10s1}) can be expressed as
\begin{align*}
\sum_{k,m_2, m_3=0}^{\infty}\sum_{m_1=0}^{k}\frac{\wedge(m_1, m_2, m_3)}{m_1! m_2! m_3! k!}\frac{(c_{i})_{k}(-k)_{m_1}}{(c_{i})_{m_1}}(\frac{1+t}{t} x_1)^{m_1} x_{2}^{m_2} x_{3}^{m_3}(-t)^{k},
\end{align*}
Replacing $k=m_1+l$, changing the summation order and simplifying, we get
\begin{align*}
\sum_{m_1, m_2, m_3=0}^{\infty}\frac{\wedge(m_1, m_2, m_3)}{m_{1}! m_{2}! m_{3}!} x_{1}^{m_{1}}x_{2}^{m_{2}} x_{3}^{m_{3}}(1+t)^{m_1} \, _{1}F_{0}\left[^{c_{i}+m_{1}}_{\,-\,}; -t\right].
\end{align*}
Evaluating the inner $ _{1}F_{0}$-series in the above equation by binomial theorem (\ref{N}) 
\begin{align*}
_{1}F_{0}\left[^{c_{i}+m_{1}}_{\,-\,}; -t\right]= (1+t)^{-c_{i}-m_{1}}
\end{align*}
and simplifying, we get the right side of this theorem. This completes the proof.

Recalling the binomial theorem (\ref{N}), we establish another infinite summation formula of Srivastava's general triple hypergeometric function.
\begin{theorem}
The following infinite summation formula of Srivastava's general triple hypergeometric function holds true:
\begin{align}
&\sum_{k=0}^{\infty}\frac{(c_{i})_{k}}{k!}\,(\frac{t+x_1}{x_1-1})^{k}\,F^{(3)}\left [^{\,(a)\,::\,(b)\,;\,(b')\,; \,(b''): \,-k, (c^{i})\,; (c')\,; \,(c'') \,;}_{\,(e)\,::\,(g)\,; \,(g')\,;\, (g'')\,:\,\,\,(h)\,\,\,\,\,;\,(h');\,(h'')\,;} \frac{1+t}{t+x_1}  x_1, x_2, x_3\right ]\notag\\
&=\, (\frac{1-x}{t+1})^{c_{i}}\,F^{(3)}\left [ x_1, x_2, x_3\right],\label{10s2}
\end{align}
where $i=1,\dots,C$.
\end{theorem}
{\bf Proof:} The proof of this theorem is similar to Theorem \ref{t9}. 
We omit the details.
\section{Conclusion}
We have obtained several infinite summation formulas involving the Srivastava's general triple hypergeometric function. We remark that by specializing the parameters in $F^{(3)}[ x_1, x_2, x_3]$, we can deduce summation formulas for the forteen Lauricella functions \cite{GL} as well as three Srivastava's triple hypergeometric functions $H_{A}$,  $H_{B}$ and  $H_{C}$, \cite{S1,SK, SM}. We have listed only some particular cases leading to results of Lauricella functions and three Srivastava's triple hypergeometric functions.  The infinite summation formulas for the remaining Lauricella functions and three Srivastava's triple hypergeometric functions  can be worked out  analogously.\\

For example, specializing the parameters in (\ref{1s2}), (\ref{t1s2}) and (\ref{t2s2}) we get the  infinite summation formulas for  $F_{A}^{(3)}$ and $F_{D}^{(3)}$:
\begin{align}
&\sum_{k=0}^{\infty}\frac{(a)_{k}}{k!}\, t^{k}\, F^{\left(3\right)}_{A}\left({a}+k,{b_{1}},{b_{2}},{b_{3}};{c_{1}},{c_{2}},{c_{3}};{x_{1}},{x_{2}},{x_{3}}\right)\notag\\
&=(1-t)^{-a}\,F^{\left(3\right)}_{A}\left({a},{b_{1}},{b_{2}},{b_{3}};{c_{1}},{c_{2}},{c_{3}};\frac{x_{1}}{1-t},\frac{x_{2}}{1-t},\frac{x_{3}}{1-t}\right);
\end{align}
\begin{align}
&\sum_{k=0}^{\infty}\frac{(a)_{k}}{k!}\, t^{k}\, F^{\left(3\right)}_{D}\left({a}+k,{b_{1}},{b_{2}},{b_{3}};{c};{x_{1}},{x_{2}},{x_{3}}\right)\notag\\
&=(1-t)^{-a}\,F^{\left(3\right)}_{D}\left({a},{b_{1}},{b_{2}},{b_{3}};{c};\frac{x_{1}}{1-t},\frac{x_{2}}{1-t},\frac{x_{3}}{1-t}\right).
\end{align}
 Again, specializing the parameters in (\ref{s3}) we obtain the  infinite summation formulas for  $H_{A}$:
\begin{align}
&\sum_{k=0}^{\infty}\frac{(a)_{k}(b_{1})_{k}}{(c_{1})_{k}\,k!}\, t^{k}\,H_{A}(a+k, b_1+k, b_2; c_1+k, c_2 ; x_1,x_2,x_3)\notag\\
&=\,H_{A}(a, b_1, b_2; c_1, c_2;\,{x_{1}+t} ,x_2,{x_{3}}).
\end{align}

\end{document}